\documentclass[12pt]{article}
\usepackage[top=1in,bottom=1in,left=1in,right=1in]{geometry}
\usepackage{indentfirst}
\usepackage{latexsym}
\usepackage{graphicx}
\usepackage{mathrsfs}
\usepackage{amsfonts,amsmath,amsthm,amssymb,mathtools}
\usepackage{longtable, latexsym, euscript, epic, eepic, color}

\newtheorem{lemma}{Lemma}[section]
\newtheorem{theorem}{Theorem}[section]

\newtheorem{proposition}{Proposition}[section]

\newtheorem{remark}{Remark}

\newtheorem{hypothesis}{Hypothesis}

\begin{document}
\baselineskip=19pt

\title{Flag-transitive point-primitive non-symmetric $2$-$(v,k,2)$ designs with alternating socle\footnote{This work is
supported by the National Natural Science Foundation
of China (Grant No.11471123).}}

\author{Hongxue  Liang, Shenglin Zhou
%\small \it School of Mathematics, South China University of Technology,\\
%\small \it Guangzhou 510640, P.R. China
%\date{}
}
\maketitle

\begin{abstract}
This paper studies flag-transitive point-primitive non-symmetric $2$-($v,k,2$) designs.
We prove that if $\mathcal{D}$ is a non-trivial non-symmetric $2$-$(v,k,2)$ design
admitting a flag-transitive point-primitive automorphism group $G$ with $Soc(G)=A_{n}$ for $n\geq5$,
then $\mathcal{D}$ is a $2$-$(6,3,2)$ or $2$-$(10,4,2)$ design.

\medskip
\noindent{\bf MR(2000) Subject Classification:} 05B05, 05B25,
20B25


\end{abstract}

\section{Introduction}

A \emph{$2$-$(v,k,\lambda)$ design} is a finite incidence
structure $\cal D$=($\cal P$, $\cal B$)
consisting of $v$ points and $b$ blocks such that every
block is incident with $k$ points, every point is incident
with $r$ blocks, and any two distinct points are incident
with exactly $\lambda$ blocks. The design $\cal D$ is
called \emph{symmetric} if $v=b$ (or equivalently $r=k$) and \emph{non-trivial} if $1<k<v$.
A \emph{flag} of $\cal D$ is an incident point-block
pair $(\alpha, B)$ where $\alpha$ is a point and $B$ is a block.
An automorphism of $\cal D$ is a
permutation of the points which also permutes the blocks. The set
of all automorphisms of $\cal D$ with the composition of maps is a
group, denoted by $Aut({\cal D})$. A subgroup $G\leq Aut({\cal D})$
is called \emph{point-primitive} if it acts primitively on $\cal P$
and \emph{flag-transitive} if it acts transitively on
the set of flags of $\cal D$.

 It is shown in \cite{Re05} that  the socle of the automorphism group of  a flag-transitive, point-primitive symmetric $2$-$(v,k,2)$ design cannot be alternating or sporadic.
 Motivated by this article, it is natural to consider the case of non-symmetric designs. Recently, we proved in \cite{LiangZhou} that, for a non-symmetric $2$-($v,k,2$) design, if $G\leq Aut(\cal D)$ is flag-transitive and point-primitive then $G$ must be an affine or almost simple group. Moreover, if the socle of $G$ is sporadic, then $\cal D$ is the unique $2$-$(176,8,2)$ design with $G=HS$, the Higman-Sims simple group.
This article is a continuation of \cite{LiangZhou}. Here we treat almost simple groups in which $Soc(G)$ is an alternating group.
The main result is the following.

\begin{theorem}\label{th1}
If $\mathcal{D}$ is a non-trivial non-symmetric $2$-{\rm(}$v,k,2${\rm)} design
admitting a flag-transitive point-primitive automorphism group $G$ with alternating socle $A_n$ for $n\geq 5$, then
%then up to isomorphism $(\mathcal{D},G)$ is one of the following :
\begin{enumerate}
\item[\rm(i)] \, $\mathcal{D}$ is a unique $2$-$(6,3,2)$ $design$ and $G=A_{5}$, or

\item[\rm(ii)] \, $\mathcal{D}$ is a unique $2$-$(10,4,2)$ $design$ and $G=A_{5}$, $A_{6}$ or $S_{6}$.
\end{enumerate}

\end{theorem}
The structure of our paper is as follows. In Section \ref{sec2},
we give some preliminary lemmas on flag-transitive designs and
permutation groups.
In Section \ref{sec3}, we complete the proof of Theorem \ref{th1} in 5 parts.

\section{Preliminaries}\label{sec2}
\begin{lemma} \label{condition 1}
The parameters $v, b, r, k, \lambda$ of a non-trivial $2$-$(v,k,\lambda)$ design
satisfy the following arithmetic conditions:
\begin{enumerate}
\item[\rm(i)] \, $vr=bk$;

\item[\rm(ii)] \, $\lambda(v-1)=r(k-1)$;

\item[\rm(iii)] \, $b\geq v$ and $k\leq r$.
\end{enumerate}
\end{lemma}
\begin{flushleft}
In particular, if the design is non-symmetric then $b>v$ and $k<r$.
\end{flushleft}

\begin{lemma} \label{condition 2}
Let $\cal D$ be a non-trivial $2$-$(v,k,\lambda)$ design.
Let $\alpha$ be a point of $\cal D$ and $G$ be a flag-transitive
automorphism group of $\cal D$.
\begin{enumerate}
\item[\rm(i)] \, $r^{2}>\lambda v$
                 and $|G_\alpha|^3>\lambda |G|$. In particular, $r^2>v$.
\item[\rm(ii)] \, $r\mid\lambda(v-1,|G_{\alpha}|)$, where
                  $G_{\alpha}$ is the stabilizer of $\alpha$.
\item[\rm(iii)] \, If $d$ is any non-trivial subdegree of $G$, then
                   $r\mid\lambda d$ {\rm (}and so $\frac{r}{(r,\lambda)}\mid d${\rm )}.
\end{enumerate}
\end{lemma}

{\bf Proof}.
(i) The equality $r=\frac{\lambda(v-1)}{k-1}$ implies
$\lambda v=r(k-1)+\lambda\leq r(r-1)+\lambda=r^2-r+\lambda$, and
the non-triviality of $\cal D$ implies $r>\lambda$,
and so $r^2>\lambda v$. Combining this with
$v=|G:G_\alpha|$ and $r\leq|G_\alpha|$ by the flag-transitivity
of $G$, we have $|G_\alpha|^3>\lambda|G|$.
(ii) Since $G$ is flag-transitive and $\lambda(v-1)=r(k-1)$,
we have $r\mid\lambda(v-1)$ and $r\mid|G_\alpha|$.
It follows that $r$ divides $(\lambda(v-1),|G_\alpha|)$,
and hence $r\mid\lambda(v-1,|G_\alpha|)$.
Part (iii) was proved in \cite[p.91]{Dav1} and \cite{Dav2}.

$\hfill\square$

\begin{lemma}{\rm (\cite[p.366]{Lie1})}\label{3 cases}
If $G$ is $A_{n}$ or $S_{n}$, acting on a set $\Omega$ of size $n$, and
$H$ is any maximal subgroup of $G$ with $H\neq A_{n}$, then $H$ satisfies one of the following:
\begin{enumerate}
\item[\rm(i)] \, $H=(S_{\ell}\times S_{m})\cap G$, with $n=\ell+m$ and $\ell\neq m$ {\rm(}intransitive case{\rm)};

\item[\rm(ii)] \, $H=(S_{\ell}\wr S_{m})\cap G$, with $n=\ell m$, $\ell>1$ and $\ell\neq m$ {\rm(}imprimitive  case{\rm)};

\item[\rm(iii)] \, $H=AGL_{m}(p)\cap G$, with $n=p^{m}$ and $p$ a prime {\rm(}affine case{\rm)};

\item[\rm(iv)] \, $H=(T^{m}.(Out~T\times S_{m}))\cap G$, with $T$ a nonabelian simple group, $m\geq2$
and $n=|T|^{m-1}$ {\rm(}diagonal case{\rm)};

\item[\rm(v)] \, $H=(S_{\ell}\wr S_{m})\cap G$, with $n=\ell^{m}$, $\ell\geq5$
and $m>1$ {\rm(}wreath case{\rm)};

\item[\rm(vi)] \, $T\unlhd H\leq Aut(T)$, with $T$ a nonabelian simple group, $T\neq A_{n}$ and $H$
acting primitively on $\Omega$ {\rm(}almost simple case{\rm)}.
\end{enumerate}
\end{lemma}

\begin{remark} {\rm This lemma does not deal with the groups $M_{10}$, $PGL_{2}(9)$
and $P\Gamma L_{2}(9)$ that have $A_{6}$ as socle. These exceptional cases will
be handled in the first part of Section \ref{sec3}.}
\end{remark}

\begin{lemma}{\rm (\cite[Theorem (b)(I)]{Lie2})} \label{odd degree}
Let G be a primitive permutation group of odd degree n on a set $\Omega$ with simple socle
$X= Soc(G)$, and let $H=G_{\alpha}$, $\alpha\in \Omega$. If $X\cong A_{c}$, an alternating group, then
one of the following holds:
\begin{enumerate}
\item[\rm(i)] \, H is intransitive, and $H=(S_{a}\times S_{c-a})\cap G$ where $1\leq a <\frac{1}{2}c$;

\item[\rm(ii)] \, H is transitive and imprimitive, and $H=(S_{a}\wr S_{c/a})\cap G$ where $a>1$
and $a \mid c$;

\item[\rm(iii)] \,H is primitive, $n=15$ and $G\cong A_{7}$.
\end{enumerate}
\end{lemma}

\begin{lemma}{\rm (\cite[Theorem 5.2A]{Dix1996})}\label{3 inequations}
Let $G=Alt(\Omega)$ where $n=|\Omega|\geq5$, and let $s$ be an integer
with $1\leq s\leq \frac{n}{2}$.
Suppose that, $K\leq G$ has index $|G:K|<\binom{n}{s}$. Then one of the following holds:
\begin{enumerate}
\item[\rm(i)] \, For some $\Delta\subset\Omega$ with $|\Delta|<s$ we have $G_{(\Delta)}\leq K \leq G_{\{\Delta\}}$;

\item[\rm(ii)] \, $n=2m$ is even, $K$ is imprimitive with two blocks of size $m$, and $|G:K|=\frac{1}{2}\binom{n}{m}$; or

\item[\rm(iii)] \, one of six exceptional cases hold where:
\begin{enumerate}
\item[\rm(a)] \, $K$ is imprimitive on $\Omega$ and $(n, s, |G:K|)=(6, 3, 15)$;

\item[\rm(b)] \, $K$ is primitive on $\Omega$ and $(n, s, |G:K|,K)=(5, 2, 6, 5:2)$,
$(6, 2, 6, PSL_{2}(5))$, $(7, 2, 15, PSL_{3}(2))$, $(8, 2, 15, AGL_{3}(2))$ or $(9, 4, 120, P\Gamma L_{2}(8))$.
\end{enumerate}
\end{enumerate}
\end{lemma}

\begin{remark}
 {\rm (1) From part (i) of Lemma \ref{3 inequations} we know that $K$ contains the
alternating group $G_{(\Delta)}=Alt(\Omega \backslash \Delta)$ of degree $n-s+1$.

(2) A result similar to Lemma \ref{3 inequations} holds for the finite symmetric groups $Sym(\Omega)$
 which can be found in \cite[Theorem 5.2B]{Dix1996}.}
 \end{remark}

We will also need some elementary inequalities.
\begin{lemma}\label{inequation}
Let $s$ and $t$ be two positive integers.
\begin{enumerate}
\item[\rm(i)] \, If $t>s\geq 7$, then $\binom{s+t}{s}>2t^{4}>2s^{2}t^{2}$.

\item[\rm(ii)] \, If $s\geq 6$ and $t\geq2$, then $2^{(s-1)(t-1)}>2s^{4}\binom{t}{2}^{2}$
implies $2^{s(t-1)}>2(s+1)^{4}\binom{t}{2}^{2}$.

\item[\rm(iii)] \, If $t\geq 6$ and $s\geq2$, then $2^{(s-1)(t-1)}>2s^{4}\binom{t}{2}^{2}$
implies $2^{(s-1)t}>2s^{4}\binom{t+1}{2}^{2}$.

\item[\rm(iv)] \, If $t\geq4$, and $s\geq3$, then $\binom{s+t}{s}>2s^{2}t^{2}$ implies
$\binom{s+t+1}{s}>2s^{2}(t+1)^{2}$.
\end{enumerate}
\end{lemma}

{\bf Proof.}\, (i) It is necessary to prove that $\binom{s+t}{s}>2t^{4}$ holds.
Since $t>s\geq 7$ then $[\frac{s+t}{2}]\geq s \geq 7$, it follows that
$\binom{s+t}{s}\geq\binom{t+7}{7}>2t^{4}$.

(ii) Suppose that $s\geq 6$, $t\geq2$ and $2^{(s-1)(t-1)}>2s^{4}\binom{t}{2}^{2}$.
Then $$2^{s(t-1)}=2^{(s-1)(t-1)}2^{t-1}>2s^{4}\binom{t}{2}^{2}2^{t-1}
=2(s+1)^{4}\binom{t}{2}^{2}(1-\frac{1}{s+1})^{4}2^{t-1}.$$
Combing this with the fact $(1-\frac{1}{s+1})^{4}2^{t-1}\geq2\times(\frac{6}{7})^{4}>1$ gives (ii).

(iii) Suppose that $t\geq 6$, $s\geq2$ and $2^{(s-1)(t-1)}>2s^{4}\binom{t}{2}^{2}$.
Then $$2^{(s-1)t}=2^{(s-1)(t-1)}2^{s-1}>2s^{4}\binom{t}{2}^{2}2^{s-1}
=2s^{4}\binom{t+1}{2}^{2}(1-\frac{2}{t+1})^{2}2^{s-1}.$$
Combing this with the fact $(1-\frac{2}{t+1})^{2}2^{s-1}\geq2\times(\frac{5}{7})^{2}>1$ gives (iii).

(iv) Suppose that $t\geq 4$, $s\geq3$ and $\binom{s+t}{s}>2s^{2}{t}^{2}$. Then
$$\binom{s+t+1}{s}=\binom{s+t}{s}\frac{s+t+1}{t+1}>2s^{2}{t}^{2}\frac{s+t+1}{t+1}
=2s^{2}(t+1)^{2}\frac{(s+t+1)t^{2}}{(t+1)^{3}}.$$
The fact that $(s+t+1)t^{2}>(t+1)^{3}$ gives (iv).

$\hfill\square$

\section{Proof of Theorem \ref{th1}}\label{sec3}

In this section, $\mathcal{D}$ denotes a non-trivial non-symmetric 2-$(v,k,2)$
design if without special statement, and $G\leq Aut({\cal D})$ is
flag-transitive point-primitive with $Soc(G)=A_{n}$.
Let $\alpha$ be a point of $\mathcal{D}$ and $H=G_{\alpha}$.
Since $G$ is point-primitive, $H$ is a maximal subgroup of $G$
by \cite[Theorem 8.2]{Wie1964}.
Furthermore, by the flag-transitivity of $G$, we have that $v=|G:H|$,
$b\mid |G|$, $r\mid |H|$ and $r^{2}>2v$ by Lemma \ref{condition 2} (i).

If $r$ is odd, then $(r,2)=1$. This case has been classified by Zhou and Wang in \cite{yajie2015}, and we get the following:

\begin{proposition} \label{pro1}
Let $\mathcal{D}$ be a non-trivial non-symmetric $2$-$(v,k,2)$ design
admitting a flag-transitive point-primitive automorphism
group $G$ with $Soc(G)=A_{n}$, $n\geq 5$. If the replication number $r$ is odd, then $\cal D$ is the unique
$2$-$(6,3,2)$ design  and $G=A_{5}$.
\end{proposition}

%\subsection{$r$ is even}\label{sub2}
Now we assume that $r$ is even in the following.

Suppose first that $n=6$ and $G\cong M_{10}$, $PGL_{2}(9)$ or $P\Gamma L_{2}(9)$.
Each of these groups has exactly three maximal subgroups with index greater than 2,
and their indices are precisely 45, 36 and 10.
By using the computer algebra system {\sf GAP} \cite{gap},
for $v=45$, 36 or 10, we will compute the parameters $(v,b,r,k)$ that satisfy the following conditions:
\begin{equation}
r\mid (2(v-1),|H|);
\end{equation}
\begin{equation}
   r^{2}>2v;
\end{equation}
\begin{equation}
   2\mid r;
\end{equation}
\begin{equation}
   r(k-1)=2(v-1);
\end{equation}
\begin{equation}
    r>k>2;
\end{equation}
\begin{equation}
   b=\frac{vr}{k};
\end{equation}
We obtain two possible parameters $(v,b,r,k)$ as follows:
\begin{center}
  $(10, 15, 6, 4)$ and $(36, 45, 10, 8)$.
\end{center}
Now we consider the existence of flag-transitive point-primitive non-symmetric designs
with above possible parameters.

Suppose first that there exists a $2$-$(10,4,2)$ design $\mathcal{D}$ with a flag-transitive
point-primitive automorphism group $G$.
Let the point set $\mathcal{P}=\{1, 2, ..., 10\}$,
and the group $G=M_{10}$, $PGL_{2}(9)$ or $P\Gamma L_{2}(9)$ as the primitive permutation
group of degree 10 acting on $\mathcal{P}$.
Since $G$ is flag-transitive, $G$ acts block-transitively on $\mathcal{B}$,
so $|G|/b=|G_{B}|$, where $B$ is a block. For each case, using the
command {\tt Subgroups(G:OrderEqual:=n)}
where ${\tt n}=|G|/b$ by {\sf Magma} \cite{Magma},
it turns out that
$G$ has no subgroup of order ${\tt n}$,
which contradicts the fact that $G_{B}$ is a subgroup of order $|G|/b$.

Assume next that there exists a $2$-$(36,8,2)$ design $\mathcal{D}$ with the flag-transitive
point-primitive automorphism group $G=M_{10}$, $PGL_{2}(9)$ or $P\Gamma L_{2}(9)$.

When $(v,G)=(36,P\Gamma L_{2}(9))$, by the {\sf Magma}-command {\tt Subgroups(G:OrderEqual:=n)}
where ${\tt n}=|G|/b$, we get the block stabilizer $G_{B}$. Since $G$ is flag-transitive,
$G_{B}$ is transitive on $B$, and so $B$ is an orbit of $G_{B}$ acting on the point set $\mathcal{P}$.
Using the {\sf Magma}-command {\tt Orbits(GB)} where ${\tt GB}=G_{B}$,
it turns out that $G_{B}$ has no orbit of length $k$, a contradiction.

Now assume that $(v,G)=(36,M_{10})$ or $(36,PGL_{2}(9))$.
By the definition of $2$-$(v,k,\lambda)$ design, every pair of distinct points is contained
in exactly $\lambda$ blocks for some positive constant $\lambda$, i.e. \emph{pairwise balanced}.
However, for each case,
the command {\tt PairwiseBalancedLambda(D)} turns out false.

Take $(v,G)=(36,M_{10})$ for example, the orbits of $G_{B}$ are:
\begin{align*}
  \Delta_{0}=&\{3, 17, 18, 21\}, \\
  \Delta_{1}=&\{1, 4, 12, 14, 16, 22, 26, 34\},\\
  \Delta_{2}=&\{2, 6, 7, 9, 15, 23, 29, 36\},\\
  \Delta_{3}=&\{5, 8, 10, 11, 13, 19, 20, 24, 25, 27, 28, 30, 31, 32, 33, 35\}.
\end{align*}

As $k=8$, we take the orbit of length 8 as the block $B$,
i.e. $B=\Delta_{1}$ or $B=\Delta_{2}$.
Using the {\sf GAP}-command {\tt D1:=BlockDesign(36,[[1,4,12,14,16,22,26,34]],G)},
we obtain that $|\Delta_{1}^G|=45=b$. We take $\mathcal{P}=\{1, 2, \ldots, 36\}$,
$B=\Delta_{1}$ and $\mathcal{B}=B^G$.
Now, we just need check that each pair of distinct points is contained in 2 blocks.
However, {\tt PairwiseBlancedLambda(D1)} shows that this is not true, and so $B\neq\Delta_{1}$.
Similarly, we can get $B\neq\Delta_{2}$.
So the case $(v,G)=(36,M_{10})$ cannot occur.

Now we consider $G=A_{n}$ or $S_{n}$ with $n\geq5$.
The point stabilizer $H=G_{\alpha}$
acts both on $\mathcal{P}$ and the set $\Omega_{n}=\{1,2,\ldots,n\}$.
Then by Lemma \ref{3 cases} one of the following holds:
\begin{enumerate}
\item[\rm(i)] \, $H$ is primitive in its action on $\Omega_{n}$;

\item[\rm(ii)] \, $H$ is transitive and imprimitive in its action on $\Omega_{n}$;

\item[\rm(iii)] \, $H$ is intransitive in its action on $\Omega_{n}$.
\end{enumerate}

We analyse each of these actions separately. First of all, we assume a hypothesis.

\begin{hypothesis}\label{hypo 1}
Let $\mathcal{D}$ be a non-trivial non-symmetric $2$-$(v,k,2)$ design
admitting a flag-transitive point-primitive automorphism
group $G$ with $Soc(G)=A_{n}$ {\rm(}$n\geq 5${\rm)} and let the replication number $r$ be even.
\end{hypothesis}

\subsection{$H$ acts primitively on $\Omega_{n}$}\label{sub1}

\begin{proposition} \label{pro1}
Assume that Hypothesis \ref{hypo 1} holds and the point stabilizer $H$ acts primitively
on $\Omega_{n}$, then there exist 10 possible parameters $(n,v,b,r,k)$ which are
listed in Table \ref{tab3}.
\end{proposition}

{\bf Proof.}\,
We claim that $2\| r$. Otherwise $4\mid r$,
from the basic equation $r(k-1)=2(v-1)$, we have that $2\mid (v-1)$,
and so $v$ is odd. Thus by Lemma \ref{odd degree}, $v=15$, $G=A_{7}$ and $|H|=|G|/v=168$.
Since $r\mid (2(v-1),|H|)$, $r^{2}>2v$ and $k\geq3$, $r=7$ or $14$, which contradicts $4\mid r$.

Thus $2\| r$. Since $r>2$, there exists an odd prime $p$ that divides $r$, then $p\mid (v-1)$, and so $(p,v)=1$.
Thus $H$ contains a Sylow $p$-subgroup $P$ of $G$. Let $g\in G$ be a $p$-cycle, then there is a conjugate
of $g$ belongs to $H$. This implies that $H$ acting on $\Omega_{n}$ contains an even permutation with exactly one cycle of length $p$ and $n-p$ fixed points. By a result of Jordan \cite[Theorem 13.9]{Wie1964},
we have that $n-p\leq2$. Therefore, $n-2\leq p\leq n$, $p^{2}\nmid |G|$, and so $p^{2}\nmid r$.
It follows that $r=2(n-2)$, $2(n-1)$, $2n$ or $2n(n-2)$, where $n$, $n-1$ and $n-2$ are odd primes. Moreover, the primitivity of
$H$ acting on $\Omega_{n}$ and $H\ngeqslant A_{n}$ imply that
$v\geq \frac{[\frac{n+1}{2}]!}{2}$ by \cite[Theorem 14.2]{Wie1964}.
Combining with $r^{2}>2v$,
we have that $$r^{2}>[\frac{n+1}{2}]!.$$
Therefore, $(n,r)=$ (5, 6), (5, 10), (5, 30), (6, 10), (7, 10), (7, 14), (7, 70),
(8, 14), (9, 14) or (13, 286).
By Lemmas \ref{condition 1} and \ref{condition 2}, the facts $v\geq\frac{[\frac{n+1}{2}]!}{2}$
and $[b,v]\mid|G|$, we obtain 10 possible parameters $(n,v,b,r,k)$:
\begin{align*}
&(5, 10, 15, 6, 4),~(6, 16, 40, 10, 4),~(6, 36, 45, 10, 8),~(7, 15, 70, 14, 3),~(7, 16, 40, 10, 4),  \\
& (7, 21, 42, 10, 5),~(7, 36, 45, 10, 8),~ (7, 36, 84, 14, 6),~ (8, 15, 70, 14, 3),~(8, 36, 84, 14, 6).
\end{align*}
And we list them in Table \ref{tab3}.

$\hfill\square$

\subsection{$H$ acts transitively and imprimitively on $\Omega_{n}$}\label{sub2}

\begin{proposition} \label{pro2}
Assume that Hypothesis \ref{hypo 1} holds and the point stabilizer
$H$ acts transitively but imprimitively on $\Omega_{n}$,
then there exist 2 possible parameters $(n, v, b, r, k)=(6, 10, 15, 6, 4)$ or $(10, 126, 1050, 50, 6)$ which are listed in Table \ref{tab3}.
\end{proposition}

{\bf Proof.}\,
Suppose on the contrary that $\Sigma=\{\Delta_{0},\Delta_{1},\ldots,\Delta_{t-1}\}$
is a non-trivial partition of $\Omega_{n}$ preserved by $H$, where $|\Delta_{i}|=s$,
$0\leq i\leq t-1$, $s,t\geq2$ and $st=n$.
Then $$v=\binom{ts-1}{s-1}\binom{(t-1)s-1}{s-1}\ldots \binom{3s-1}{s-1}\binom{2s-1}{s-1}.$$

Moreover, the set $O_{j}$ of $j$-cyclic partitions with respect to $X$
(a partition of $\Omega_{n}$ into $t$ classes each of size $s$) is an union
of orbits of $H$ on $\mathcal{P}$ for $j=2,\ldots,t$ (see \cite {A.D, H.D} for definitions and details).

{\bf Case (1):}~ Suppose first that $s=2$, then $t\geq3$, $v=(2t-1)(2t-3)\cdots 5\cdot3$, and
$$d_{j}=|O_{j}|=\frac{1}{2}\binom{t}{j}\binom{s}{1}^{j}=2^{j-1}\binom{t}{j}.$$
We claim that $t<7$. If $t\geq7$, then
$v=(2t-1)(2t-3)\cdots 5\cdot 3>5t^{2}(t-1)^{2}$.
On the other hand, since $r$ divides $2d_{2}=2t(t-1)$, $2t(t-1)\geq r$,
and so $v<2t^{2}(t-1)^{2}$, a contradiction.
Thus $t<7$. For $t=3$, 4, 5 or 6, we calculate $d=2{\rm gcd}(d_{2}, d_{3})$,
which are listed in Table \ref{tab1} below.
\begin{table}[h]
\begin{center}
\caption{Possible $d$ when $s=2$.}\label{tab1}
\vspace{3mm}
\begin{tabular}{cccccc}
\hline
 $t$ & $n$ & $v$ &  $d_{2}$ &  $d_{3}$&  $d$ \\
\hline
 $3$ & $6$ & $15$   &$6$  &$4$  &$4$\\
 $4$ & $8$ & $105$  &$12$ &$16$ &$8$ \\
 $5$ & $10$ &$945$  &$20$ &$40$ &$40$\\
 $6$ & $12$ &$10395$&$30$ &$80$ &$20$\\
\hline
\end{tabular}
\end{center}
\end{table}

In each line $r\leq d$ which contradicts the fact $r^2>2v$.

{\bf Case (2):}~ Thus $s\geq3$. So $O_{j}$ is an orbit of $H$ on $\mathcal{P}$,
and $d_{j}=|O_{j}|=\binom{t}{j}\binom{s}{1}^{j}=s^{j}\binom{t}{j}$.
In particular, $d_{2}=s^{2}\binom{t}{2}$
and $r\mid 2d_{2}$. Moreover, from
$\binom{is-1}{s-1}=\frac{is-1}{s-1}\cdot\frac{is-2}{s-2}\cdots\frac{is-(s-1)}{1}>i^{s-1}$,
for $i=2,3,\ldots,t$, we have that $v>2^{(s-1)(t-1)}$.
Then
$$2\cdot2^{(s-1)(t-1)}< 2v < r^{2} \leq4s^{4}\binom{t}{2}^{2},$$
and so
$$2^{(s-1)(t-1)}<2s^{4}\binom{t}{2}^{2}.$$

We will calculate all pairs $(s,t)$ satisfying above inequality.
Since $2^{(6-1)(6-1)}=2^{25}>2\cdot6^{4}\cdot\binom{6}{2}^{2}=2^{5}\cdot3^{6}\cdot5^{2}$,
the pair $(s,t)=(6,6)$ does not satisfy the inequality, but satisfies the
conditions of Lemma \ref{inequation} (ii) and (iii).
Thus, we have either $s<6$ or $t<6$. It is not hard to get 36 pairs $(s,t)$ satisfying
the inequality:
\begin{align*}
&(3,2),~(3,3),~(3,4),~(3,5),~(3,6),~(3,7),~(3,8),~(3,9),~(3,10),~(4,2),~(4,3),~(4,4),~(4,5),\\
&(4,6),~(5,2),~(5,3),~(5,4),~(5,5),~(6,2),~(6,3),~(6,4),~(7,2),~(7,3),~(8,2),~(8,3),~(9,2),\\
&(9,3),~(10,2),~(11,2),~(12,2),~(13,2),~(14,2),~(15,2),~(16,2),~(17,2),~(18,2).
\end {align*}

For each pair $(s,t)$, we calculate the parameters $(v,b,r,k)$ satisfying Lemmas \ref{condition 1}, \ref{condition 2}, $2 \mid r$ and $r \mid2d_{2}$.
There exist 2 possible parameters $(n,v,b,r,k)$ corresponding to $(s,t)$:
\begin{align*}
&(s,t)=(3,2)~~{\rm with}~~ (n,v,b,r,k)=(6,10,15,6,4),\\
&(s,t)=(5,2)~~{\rm with}~~ (n,v,b,r,k)=(10,126,1050,50,6),
\end{align*}
which are listed in Table \ref{tab3}.

$\hfill\square$

\subsection{$H$ acts intransitively on $\Omega_{n}$}\label{sub3}

\begin{proposition}\label{pro3}
Assume that Hypothesis \ref{hypo 1} holds and the point stabilizer
$H$ acts intransitively on $\Omega_{n}$,
then there exist 15 possible parameters $(n,v,b,r,k)$ which are
listed in Table \ref{tab3}.
\end{proposition}

{\bf Proof.}\,
Since $H$ acts intransitively on $\Omega_{n}$,
we have $H=(Sym(S)\times Sym(\Omega_{n} \backslash S))\cap G$,
and without loss of generality, we may assume that $|S|=s<\frac{n}{2}$ by Lemma \ref{3 cases} (i).
By the flag-transitivity of $G$, $H$ is transitive on the blocks through $\alpha$, and so $H$
fixes exactly one point in $\mathcal{P}$. Since $H$ stabilizes only one $s$-subset of $\Omega_{n}$,
we can identify the point $\alpha$ with $S$. As the orbit of $S$ under $G$ consists of all the $s$-subsets
of $\Omega_{n}$, we can identify $\mathcal{P}$ with the set of $s$-subsets of $\Omega_{n}$.
So $v=\binom{n}{s}$, $G$ has rank $s+1$ and the subdegrees are:
$$d_{0}=1,~d_{i+1}=\binom{s}{i}\binom{n-s}{s-i},~i=0, 1, 2, \ldots, s-1.$$

We claim that $s\leq 6$. It follows from $r\mid 2d_{s}$ and $d_{s}=s(n-s)$
that $r\mid 2s(n-s)$. Combining this with $r^{2}>2v$,
we have that $2s^{2}(n-s)^{2}>\binom{n}{s}$. Since $s<\frac{n}{2}$
equals to $s<t=n-s$, we have $$2s^{2}t^{2}>\binom{s+t}{s}.$$

Combining this with Lemma \ref{inequation} (i), we have $s\leq6$.

{\bf Case (1):} ~If $s=1$, then $v=n\geq5$ and the subdegrees are 1, $n-1$.
If $k=v-1$, then $r(v-2)=2(v-1)$, and so $v-2\mid v-1$ for $(r,2)=2$, a contradiction.
Therefore, $2<k\leq v-2$. Since $G$ is $(v-2)$-transitive on $\mathcal{P}$,
$G$ acts $k$-transitively on $\mathcal{P}$, and so $b=|\mathcal{B}|=|B^{G}|=\binom{n}{k}$
for every block $B\in \mathcal{B}$. From the equality $bk=vr$, we obtain $\binom{n}{k}k=nr$.
On one hand, by $r(k-1)=2(n-1)$ and $k>2$, we have $r\leq n-1$, and so
$\binom{n}{k}k\leq n(n-1)$; on the other hand, by $2<k\leq n-2$, we have $n-i\geq k-i+2>k-i+1$
for $i=2,3,\ldots,k-1$. Thus,
$$\binom{n}{k}k=n(n-1)\cdot\frac{n-2}{k-1}\cdot\frac{n-3}{k-2}\cdots\frac{n-k+1}{2}>n(n-1),$$
a contradiction.

{\bf Case (2):} ~If $s=2$, then $v=\frac{n(n-1)}{2}$ and the subdegrees are $1$, $\binom{n-2}{2}$, $2(n-2)$.
By Lemma \ref{condition 2} (iii), $r\mid 2(\binom{n-2}{2},2(n-2))$.

$(a)$ If $n\equiv0$(mod 4) or $n\equiv2$(mod 4), then  $r\mid 2(\binom{n-2}{2},2(n-2))=n-2$,
it follows that $n(n-1)=2v<r^{2}\leq(n-2)^{2}$, which is impossible.

$(b)$ If $n\equiv1$(mod 4), then $r\mid 2(\binom{n-2}{2},2(n-2))=2(n-2)$.

Let $r=\frac{2(n-2)}{u}$ for some integer $u$.
Since $r^{2}>2v$, we have $4>\frac{4(n-2)^{2}}{n(n-1)}>u^{2}$,
which forces $u=1$. Therefore, $r=2(n-2)$. By Lemma \ref{condition 1}, $k=\frac{n+3}{2}$
and $b=\frac{2n(n-1)(n-2)}{n+3}$. Since $k$ and $b$ are integers, $(n+3)\mid 120$ with $n$ odd.
Combining this with $n\equiv1$(mod 4),
we get that $n$=5, 9, 17, 21, 37, 57 or 117.
For each $n$, we calculate the parameters $(v,b,r,k)$.
We find that if $n\in \{17, 21, 37, 57, 117\}$,
then $|G:G_B|=b<\binom{n}{3}$.
By Lemma \ref{3 inequations} and \cite[Theorem 5.2B]{Dix1996},
it is easy to know that
G has no subgroup of index $b$, a contradiction.
So we obtain 2 possible parameters $(n,v,b,r,k)$:
$$(5,10,15,6,4),~(9,36,84,14,6).$$

$(c)$ If $n\equiv3$(mod 4), then $r\mid 2(\binom{n-2}{2},2(n-2))=4(n-2)$.

Let $r=\frac{4(n-2)}{u}$ for some integer $u$.
Since $r^{2}>2v$, we have $16>\frac{16(n-2)^{2}}{n(n-1)}>u^{2}$,
and so $u=1$, 2 or 3.

If $u=1$, then $r=4(n-2)$, $k=\frac{n+5}{4}$ and $b=\frac{8n(n-1)(n-2)}{n+5}$.
As $b$ must be an integer, $(n+5)\mid 1680$. Combining with $n\equiv3$(mod 4), we
have that $n$=7, 11, 15, 19, 23, 35, 43, 51, 55, 75, 79, 107, 115, 135,
163, 235, 275, 331, 415, 555, 835 or 1675.
Similarly, by Lemma \ref{3 inequations} and \cite[Theorem 5.2B]{Dix1996},
$n\in \{7,11,15,19,23,35,43\}$
and we obtain 7 possible parameters $(n,v,b,r,k)$:
\begin{align*}
&(7,21,140,20,3),~  (11,55,495,36,4),~(15,105,1092,52,5),~ (19,171,1938,68,6),\\
&(23,253,3036,84,7),~ (35,595,7854,132,10),~ (43,903,12341,164,12).
\end{align*}

If $u=2$, then $r=2(n-2)$, $k=\frac{n+3}{2}$ and $b=\frac{2n(n-1)(n-2)}{n+3}$, and so $(n+3)\mid 120$. Combining with $n\equiv3 $(mod 4), we have that $n=$ 7 or 27.
By Lemma \ref{3 inequations} and \cite[Theorem 5.2B]{Dix1996}, $n\neq 27$,
so we obtain a possible parameter $(n,v,b,r,k)=(7,21,42,10,5)$.

If $u=3$, then $r=\frac{4(n-2)}{3}$, $k=\frac{3n+7}{4}$ and $b=\frac{8n(n-1)(n-2)}{3(3n+7)}$, and so $(3n+7)\mid 7280$.
 The facts that $r$ is an integer and $n\equiv3 $(mod 4) imply $n=$ 11, 35, 119 or 1211.
For each $n$,
we find $|G:G_B|=b<\binom{n}{3}$.
By Lemma \ref{3 inequations} and \cite[Theorem 5.2B]{Dix1996},
it is easy to know that
$G$ has no subgroup of index $b$, a contradiction.

{\bf Case (3):} ~Suppose that $3\leq s \leq6$. Now for each value of $s$, by the inequality
$2s^{2}t^{2}>\binom{s+t}{s}$ and Lemma \ref{inequation} (iv), we know that $t$ (and hence $n$) is bounded.
For example, let $s=3$, since $\binom{3+102}{3}>2\cdot3^{2}\cdot102^{2}$, we must have $4\leq t \leq101$,
and so $7\leq n \leq104$. The bounds of $n$ are listed in Table \ref{tab2} below.

\begin{table}[h]
\begin{center}
\caption{Bounds of $n$ when $3\leq s \leq6$}\label{tab2}
\vspace{3mm}
\begin{tabular}{ccc}
\hline
$s $ & $t$ & $n$ \\
\hline
$3$  & $4\leq t \leq101$ & $7\leq n \leq104$ \\
$4$  & $5\leq t \leq22$ & $9\leq n \leq26$  \\
$5$  & $6\leq t \leq12$ & $11\leq n \leq17$ \\
$6$  & $7, 8, 9$ & $13, 14, 15$ \\
\hline
\end{tabular}
\end{center}
\end{table}
Note that $v=\binom{n}{s}$, and $d_{1}=\binom{n-s}{s}$, $d_{2}=s\binom{n-s}{s-1}$, $d_{3}=\binom{s}{2}\binom{n-s}{s-2}$ are three non-trivial subdegrees of $G$ acting on $\mathcal{P}$.
Therefore, the 5-tuple $(n,v,b,r,k)$ satisfies the arithmetical conditions:
(3.1)-(3.6) and $r\mid 2d_{i}$, $i\in \{1,2,3\}$.

If $s=3$, using \texttt{GAP}, it outputs five 5-tuples:
\begin{align*}
  &(13, 286, 429, 30, 20),~(14, 364, 2002, 66, 12),~(22, 1540, 6270, 114, 28),\\
  &(32, 4960, 14880, 174, 58),~(50, 19600, 39480, 282, 140).
\end{align*}

If $s=$ 4, 5 or 6, using \texttt{GAP}, there is no parameter $(n,v,b,r,k)$ satisfying
these conditions.

Thus, we totally obtain 15 possible parameters $(n,v,b,r,k)$ and list them in
Table \ref{tab3}.

$\hfill\square$
\begin{table}[h]
\begin{center}
\caption {Potential parameters}\label{tab3}
\vspace{5mm}
\begin{tabular}{lllcc}
\hline
{\sc Case} & $(v,b,r,k)$   & $Soc(G)$ or $G$              &Proposition         & Step/Reference\\
\hline
$1$ &(10, 15, 6, 4)       &$A_{5}$       & 3.2                &(ii)\\

$2$ &                     &$A_{6}$         & 3.3                &${\mathcal D}$\\

$3$ &                     &$G=A_{5}$             & 3.4                &(v) \\

$4$ &                     &$G=S_{5}$             & 3.4                &${\mathcal D}$\\

$5$ &(15, 70, 14, 3)      &$G=A_{7}$ or $A_{8}$       & 3.2                &(v) \\

$6$ &                     &$G=S_{7}$ or $S_{8}$       & 3.2                &(i) \\

$7$ &(16, 40, 10, 4)      &$A_{6},A_{7}$   & 3.2        &(i) \\

$8$ &(21, 42, 10, 5)      &$A_{7}$       & 3.2                &(ii) \\

$9$ &                     & $A_{7}$      & 3.4                &(vi)\\

$10$&(21, 140, 20, 3)     & $A_{7}$      & 3.4                &(vi)\\

$11$&(36, 45, 10, 8)      & $A_{6},A_{7}$    &3.2       &(i) \\

$12$&(36, 84, 14, 6)      & $A_{7},A_{8}$    &3.2       &(i) \\

$13$&                     & $A_{9}$                &3.4       &(iv)\\

$14$&(55, 495, 36, 4)     & $A_{11}$              &3.4       &(iv)\\

$15$&(105, 1092, 52, 5)   & $A_{15}$              &3.4       &(iii)\\

$16$&(126, 1050, 50, 6)   & $A_{10}$              &3.3       &(iii)\\

$17$&(171, 1938, 68, 6)   & $A_{19}$              &3.4       &(iv)\\

$18$&(253, 3036, 84, 7)   & $A_{23}$              &3.4       &(iii)\\

$19$&(286, 429, 30, 20)   & $A_{13}$              &3.4       &(iii)\\

$20$&(364, 2002, 66, 12)  & $A_{14}$              &3.4       &(iv)\\

$21$&(595, 7854, 132, 10) & $A_{35}$             &3.4       &(iii)\\

$22$&(903, 12341, 164, 12)& $A_{43}$              &3.4       &(iv)\\

$23$&(1540, 6270, 114, 28)& $A_{22}$              &3.4       &(iii)\\

$24$&(4960, 14880, 174, 58)& $A_{32}$             &3.4       &(vii)\\

$25$&(19600, 39480, 282, 140)& $A_{50}$           &3.4       &(iii)\\
\hline
\end{tabular}
\end{center}
\end{table}

\subsection{ Rules out  potential parameters}\label{sub4}

Now, we will rule out 23 potential cases listed in Table \ref{tab3} in several steps.

{\bf (i)}  Rules out {\sc Cases} 6, 7, 11 and 12.

The {\sf GAP}-command {\tt {PrimitiveGroup(v,nr)}} returns the primitive
group with degree $v$ in the position $nr$ in the list of the library of the primitive permutation groups.
For each {\sc Case}, the command shows that there is no primitive group corresponding to $v$.

{\bf (ii)}  Rules out {\sc Cases} 1 and 8.

Since $G$ is flag-transitive, $|H|=|G|/v$. For each case, $H$ is primitive on $\Omega_{n}$. However,
the command {\tt {PrimitiveGroup(v,nr)}}, where ${\tt v}=n$, turns out that there is no
group of order $|G|/v$.

{\bf (iii)}  Rules out {\sc Cases} 15, 16, 18, 19, 21, 23 and 25.

Since $G$ is flag-transitive, $G$ acts transitively on $\mathcal{B}$,
so $|G|/b=|G_{B}|$, where $B$ is a block. For each case, using the
{\sf Magma}-command {\tt Subgroups(G:OrderEqual:=n)} where ${\tt n}=|G|/b$, it turns out that
$G$ has no subgroup of order ${\tt n}$,
which contradicts the fact that $G_{B}$ is a subgroup of order $|G|/b$.
When ${\tt v}\geq 2500$, the {\sf GAP}-command {\tt {PrimitiveGroup(v,nr)}}
does not know the group of degree ${\tt v}$.
For {\sc Case} 25, $G=A_{50}$ or $S_{50}$,
we use the {\sf Magma}-command ${\tt G:=Alt(50)}$ or ${\tt G:=Sym(50)}$ to get the group $G$,
and {\tt Subgroups(G:OrderEqual:=n)} where ${\tt n}=|G|/b$
to get that $G$ does not have such a subgroup of order $|G|/b$, a contradiction.

{\bf (iv)}  Rules out {\sc Cases} 13, 14, 17, 20 and 22.

Since $G_{B}$ is transitive on $B$, $B$ is an orbit
of $G_{B}$ acting on the point set $\mathcal{P}$.
Using the {\sf Magma}-command ${\tt Orbits(GB)}$,
where ${\tt GB}=G_{B}$,
it turns out that $G_{B}$ has no orbit of length $b$, a contradiction.

{\bf (v)}  Rules out {\sc Cases} 3 and 5.

Using the command ${\tt Orbits(GB)}$, we get the orbits of $G_{B}$.
As $|B|=k$, we take the orbit of
length $k$ as the block $B$.
Since $G$ acts transitively on $\mathcal{B}$,
$|B^{G}|=b$. However, for each case,
using the {\sf  GAP}-command {\tt {OrbitLength(G,B,OnSets)}},
we get that $|B^{G}|<b$, a contradiction.

Take {\sc Case 3}, $(v,G)=(10,A_{5})$ for example, the orbits of $G_{B}$ are:
\begin{center}
$\Delta_0=\{1, 8\}, ~~\Delta_1=\{2, 6 \}, ~~
\Delta_2=\{3, 5\}, ~~\Delta_3=\{4, 7, 9, 10\}.$
\end{center}
So $B=\Delta_3=\{4, 7, 9, 10\}$.
However, {\tt {OrbitLength(G,[4,7,9,10],OnSets)}}
turns out $|B^{G}|=5<15$, a contradiction.

{\bf (vi)}  Rules out {\sc Cases} 9 and 10.

For each case,
the {\sf GAP}-command {\tt PairwiseBalancedLambda(D)} turns out false, which means that $D$ is not pairwise balanced, contradicting the definition of design.

For {\sc Case} 9, take $(v,G)=(21,A_{7})$ for example, the orbits of $G_{B}$ are:
$$
\begin{array}{lll}
  &\Delta_{0}=\{13\}, &\Delta_{1}=\{2, 7, 12, 14, 15\},\\
  &\Delta_{2}=\{4, 9, 16, 19, 20\}, &\Delta_{3}=\{1, 3, 5, 6, 8, 10, 11, 17, 18, 21\}.
\end{array}
$$
As $k=5$, we take the orbit of
length 5 as the block $B$, i.e.
$B=\Delta_{1}$ or $B=\Delta_{2}$.
Using the {\sf GAP}-command
${\tt D:= BlockDesign(21,[[2, 7, 12, 14, 15]],G),}$
we obtain that $|\Delta_{1}^{G}|=42$ and $\Delta_{2}\in \Delta_{1}^{G}$.
Without loss of generality, we take
$\mathcal{P}=\{1,2,\ldots,21\}$, $B=\Delta_{1}$ and $\mathcal{B}=\Delta_{1}^{G}$.
Now, we just need check that {\tt D} is pairwise balanced.
However, {\tt PairwiseBalancedLambda(D)} shows that this is not true.
So the case $(v,G)=(21,A_{7})$ cannot occur.

{\bf (vii)}  Rules out {\sc Case} 24.

Consider first that $(v,G)=(4960,A_{32})$.
Let $\Omega_{n}=\{1,2,\ldots,32\}$, then $G$ acts primitively on $\Omega_{n}$.
Let $\mathcal{P}=\Omega_{n}^{\{3\}}$ denote the set of all 3-subsets (that is, subsets of size 3) of $\Omega_{n}$.
Then $G$ acts on $\mathcal{P}$ in a natural way, namely:
$(\alpha_{1}, \alpha_{2}, \alpha_{3})^{g}=(\alpha_{1}^{g}, \alpha_{2}^{g}, \alpha_{3}^{g})$
for all $g\in G$ and $|\mathcal{P}|=\binom{32}{3}=4960$.
Using the {\sf Magma}-command ${\tt G:=Alt(32)}$ to get the group $G$, and
{\tt Subgroups(G:OrderEqual:=n)} where ${\tt n}=|G|/b$ to know that $G$ contains
only one conjugacy class of subgroups of order $|G|/b$, denoted by $K$, so the block stabilizer $G_{B}$
is conjugate to $K$, and then there is a block $B_{0}$ such that $K=G_{B_{0}}$.
Since $G$ is flag-transitive, $K$ is transitive on $B_{0}$,
that is, $B_{0}$ is an orbit of $K$ acting on $\mathcal{P}$.
Take $S=\{1,2,3\}\in \mathcal{P}$.
Using the command ${\tt OrbitLength(G,S,OnSets)}$, we have that $G$ acts transitively on $\mathcal{P}$, and using the command ${\tt OrbitLength(K,S',OnSets)}$
for all ${\tt S'}\in \mathcal{P}$, we get that $K$
acting on $\mathcal{P}$ has exactly only one orbit $\Gamma$ of length 58.
As $k=58$, we take $B_{0}=\Gamma$. Furthermore, the {\sf Magma}-command ${\tt O:=\Gamma^\wedge G}$
turns out that ${\tt|O|}=14880=b$, and then we take $\mathcal{B}={\tt O}$.
Now, we just need check that each pair of distinct points is
contained in 2 blocks. Let $S_{1}=\{1,2,3\}$, $S_{2}=\{5,6,9\}\in \mathcal{P}$.
Using {\sf Magma},
it is easy to know that there is no block in $\mathcal{B}$
containing both $S_{1}$ and $S_{2}$ ,
a contradiction. So the case $(v,G)=(4960,A_{32})$ cannot occur.

The analysis of $(v,G)=(4960,S_{32})$ is the same as above.
Let $\mathcal{P}=\Omega_{n}^{\{3\}}$ and $G$ acts on $\mathcal{P}$ in the natural way.
Using {\sf GAP} and {\sf Magma}, we get the group $G$, the subgroup $K=G_{B_{0}}$ of order $|G|/b$,
and the orbit $\Gamma$ ($|\Gamma|=58$) of $K$ acting on $\mathcal{P}$.
Since $G$ is flag-transitive, $G$ acts transitively on $\mathcal{B}$, and so $|\Gamma^{G}|=14880$.
Using the {\sf Magma}-command
\begin{align*}
{\tt M:=PermutationGroup}<32|&(1,2,3,4,5,6,7,8,9,10,11)(12,13,14,15,16,17,18,19,\\
                            &20,21,22,23,24)(25,26,27),(9,11,15),(4,13,21)>,
\end {align*}

\noindent we get a subgroup $M<G$, and so $|\Gamma^{M}|\leq14880$.
However, ${\tt O:=\Gamma^\wedge M}$ turns out that $|{\tt O}|=36432>14880$, a contradiction.
So the case $(v,G)=(4960,S_{32})$ is ruled out.

$\hfill\square$

\subsection{ The unique non-symmetric $2$-$(10,4,2)$ design}\label{sub 5}
For {\sc Case} 2, $(v,G)=(10,A_{6})$, using the {\sf GAP}-command {\tt {PrimitiveGroup(10,3)}},
we get the primitive permutation representations
of $G=A_{6}\cong PSL(2, 9)$ acting on the set $\mathcal{P}=\{1,2,\dots ,10\}$.
By {\sf Magma}, it is easy to know that $G$ contains
two conjugacy classes of subgroups of order 24, denoted by $K_{1}$ and $K_{2}$ as representatives.
As $G$ acts flag-transitively, for any $B\in\mathcal{B}$, $|G_{B}|=|G|/b=24$,
so $G_{B}$ is conjugate to either $K_{1}$ or $K_{2}$.

Assume first that $G_{B}$ is conjugate to $K_{1}$, then there exists a block $B_0$ such that $K_{1}=G_{B_{0}}$.
By the flag-transitivity of $G$,
$B_{0}$ is an orbit of $K_{1}$ acting on $\mathcal{P}$.
The {\sf Magma}-command ${\tt {Orbits(K1)}}$,
where ${\tt K1}=K_{1}$, turns out that the lengths of the orbits of $K_{1}$ are 4 and 6.
As $|B_{0}|=k=4$, we take the orbit of length 4 as the block $B_{0}$.
Using {\sf GAP}, we obtain that $|{B_{0}}^G|=15$ and every pair of distinct points
is contained in exactly 2 blocks.
Therefore, we get a point-primitive non-symmetric $2$-$(10,4,2)$ design, denoted by ${\mathcal D}_{1}$.
Since $G$ acts transitively on $\mathcal{B}$ and $G_{B_{0}}$
acts transitively on $B_{0}$, $G$ acts flag-transitively on ${\mathcal D}_{1}$.
Thus, we get a desired flag-transitive point-primitive $2$-$(10,4,2)$
non-symmetric design ${\mathcal D}_{1}$.

Consider next that $G_{B}$ is conjugate to $K_{2}$. By {\sf Magma}, we obtain that the orbits of $K_{2}$
are $\{2, 6, 7, 10\}$ and $\{1, 3, 4, 5, 8, 9\}$, and so we take $B_{0}=\{2, 6, 7, 10\}$.
Similarly, by the facts that $|B_{0}^G|=15$ and every pair of distinct points
is contained in exactly 2 blocks, we get a
flag-transitive point-primitive 2-$(10,4,2)$ non-symmetric design ${\mathcal D}_{2}$.

The analysis of {\sc Case 2}, $(v,G)=(10,S_{6})$ is the same as above, here
the group $G$ is {\tt PrimitiveGroup(10,5)}.
There are two conjugacy classes of subgroups of order 48,
denoted by $M_{1}$ and $M_{2}$ as representatives.
The orbits of $M_{1}$ are $\{2, 4, 6, 8\}$ and $\{1, 3, 5, 7, 9, 10\}$,
and so we take $B_{0}=\{2, 4, 6, 8\}$. Using {\sf GAP}, we get that $|{B_{0}}^G|=15$ and every pair of distinct points
is contained in exactly 2 blocks, and so we obtain a 2-$(10,4,2)$ design ${\cal D}_{3}$ as desired.
The lengths of the orbits of $M_{2}$ are 4 and 6, and so we take the orbit of size 4 as the block $B_{0}$,
which satisfies $|{B_{0}}^G|=15$ and $\lambda=2$, hence, we get a desired design ${\cal D}_{4}$.

For {\sc Case} 4, we get the primitive group $G=S_{5}$
by {\tt {PrimitiveGroup(10,2)}}.
Using the
command {\tt Subgroups(G:OrderEqual:=n)} where ${\tt n}=|G|/b$, we know that
there is exactly only one conjugacy class of subgroups of order 8, denoted by $L$,
so the block stabilizer $G_{B}$
is conjugate to $L$, and then there is a block $B_{0}$ such that $L=G_{B_{0}}$.
${\tt Orbits(L)}$ gives the orbits of $L$:
\begin{center}
$\Delta_0=\{3, 5\}, ~~\Delta_1=\{1, 2, 6, 8\}, ~~
\Delta_2=\{4, 7, 9, 10\}.$
\end{center}
So either $B_0=\Delta_1$ or $B_0=\Delta_2$. If $B_0=\Delta_2$,
{\tt {OrbitLength(G,[4,7,9,10],OnSets)}} shows that $|B_{0}^{G}|=5<15$, a contradiction.
If $B_0=\Delta_1$, ${\tt D:= BlockDesign(10,[[1, 2, 6, 8]],G)}$
shows that $|B_{0}^{G}|=15$, and {\tt PairwiseBalancedLambda(D)} gives $\lambda=2$.
We take $\mathcal{P}=\{1,2,\ldots,10\}$, $B_{0}=\Delta_1$ and $\mathcal{B}=B_{0}^{G}$.
Therefore, we construct a non-symmetric $2$-$(10,4,2)$ design $\mathcal{D}_{5}$,
which is flag-transitive and point-primitive.

The {\sf GAP}-command {\tt IsIsomorphicBlockDesign(D1,D2)} turns out true when ${\tt D1}={\mathcal D}_{1}$ and
${\tt D2}={\mathcal D}_{2}$, ${\mathcal D}_{3}$, ${\mathcal D}_{4}$ or ${\mathcal D}_{5}$.
So up to isomorphism we take these five designs as the same design, denoted by ${\mathcal D}$.

This completes the proof of Theorem \ref{th1}.

$\hfill\square$

\section*{Author's addresses}
\noindent{\sc Hongxue Liang}\quad 904102195@qq.com\\
{\sc Shenglin Zhou}\quad slzhou@scut.edu.cn\\
School of Mathematics\\
South China University of Technology\\
Guangzhou, Guangdong 510640\\
People's Republic of China

\end{document}